
\magnification=\magstep1
\input amstex
\documentstyle{amsppt}
\leftheadtext{E. Makai, Jr.}
\rightheadtext{Epireflective subcategories of ${\bold{TOP}}$,
${\bold{T_2UNIF}}$, ${\bold{UNIF}}$}
\topmatter
\title Epireflective subcategories of ${\bold{TOP}}$, ${\bold{T_2UNIF}}$,
${\bold{UNIF}}$,
closed under epimorphic images, or being algebraic\endtitle
\author E. Makai, Jr.* \endauthor
\address Alfr\'ed R\'enyi Mathematical Institute, 
\newline
Hungarian Academy of Sciences,
\newline
H-1364 Budapest, Pf. 127, 
\newline
HUNGARY
\newline
{\rm{http://www.renyi.mta.hu/\~{}makai}}
\endaddress
\email makai.endre\@renyi.mta.hu\endemail
\thanks *Research (partially) supported by Hungarian National Foundation for 
Scientific Research, grant nos. K68398, K75016, K81146\endthanks
\subjclass {\it Mathematics Subject Classification.} Primary: 54E
\newline
15,
Secondary: 54B30, 18C05, 18C10\endsubjclass
\abstract 
The epireflective subcategories of ${\bold{Top}}$, that are closed under
epimorphic (or bimorphic) images, are $\{ X \mid |X| \le 1 \} $, 
$\{ X \mid X$ is indiscrete$\} $ and ${\bold{Top}}$.
The epireflective subcategories of ${\bold{T_2Unif}}$, closed under
epimorphic images, are: $\{ X \mid |X| \le 1 \} $, $\{ X \mid X$ is compact
$T_2 \} $, $\{ X \mid $ covering character of $X$ is $ \le \lambda _0 \} $ 
(where $\lambda _0$ is an infinite cardinal), and ${\bold{T_2Unif}}$. 
The epireflective subcategories of ${\bold{Unif}}$, closed under
epimorphic (or bimorphic)
images, are: $\{ X \mid |X| \le 1 \} $, $\{ X \mid X$ is 
indiscrete$\} $, 
$\{ X \mid $ covering character of $X$ is $ \le \lambda _0 \} $ (where
$\lambda _0$ is an infinite cardinal), and ${\bold{Unif}}$.
The epireflective subcategories of ${\bold{Top}}$, that are
algebraic categories, are $\{ X \mid |X| \le 1 \} $, and 
$\{ X \mid X$ is indiscrete$\} $.
The subcategories of ${\bold{Unif}}$, closed under products and closed
subspaces and being varietal, are $\{ X \mid |X| \le 1 \} $, 
$\{ X \mid X$ is indiscrete$\} $,
$\{ X \mid X$ is compact $T_2 \} $.
The subcategories of ${\bold{Unif}}$, closed under products and closed
subspaces and being algebraic, are $\{ X \mid X$ is indiscrete$ \} $, and all
epireflective subcategories of $\{ X \mid X$ is compact $T_2 \} $. Also we
give a sharpened form of a theorem of Kannan-Soundararajan about classes of
$T_3$ spaces, closed for products, closed subspaces and surjective images.
\endabstract
\endtopmatter\document

\head \S 1. Preliminaries\endhead

Birkhoff's theorem in universal algebra says that varieties are
characterized, in a given type of universal algebras (i.e., 
given operations, with
given arities), as those being
closed under products, subalgebras and homomorphic
images. These properties can be investigated also in other categories,
yielding Birkhoff type theorems.

In topology it seems to have been Kannan \cite{K} 
who initiated the investigation of
simultaneously reflective and coreflective subcategories in certain
categories. If we restrict our attention to simultaneously epireflective and
monocoreflective subcategories, then under suitable hypotheses, these can be
described as those closed under products, extremal subobjects, coproducts and
extremal epi
images. This poses the question if there are theorems characterizing
subcategories of certain categories, closed under several of these operations.
Birkhoff's theorem settles one of these questions.

Herrlich \cite{H81} \S 3.2
surveyed a large number of closure operations on subcategories
of a given category, and the subcategories closed under some subsets of these
closure operations, for categories ocurring in topology. 

The category of compact $T_2$ topological spaces is characterized in several
different and nice ways, cf., e.g., de Groot's
famous characterization in \cite{W}, p. 51, Franklin-Lutzer-Thomas \cite{FLT} 
Theorem 3.10, Richter \cite{R82}, Corollary 1.7, \cite{R85a}, Corollary 4.9, 
\cite{R91b}, Theorem 4.5,
Corollary 4.7, 
\cite{R92a}, Remark 4.7.
Also cf. \cite{R91b}, Corollary 4.7, for characterizing 
epireflective subcategories of compact $T_2$ topological spaces, and 
\cite{R92a}, Remarks 2.3 and 4.7,
for characterizing reflective subcategories 
of compact $T_2$ topological spaces, containing the two-point discrete space.

The category of sets is denoted by ${\bold{Set}}$.
The categories of $T_2$ topological
spaces, $T_2$ proximity spaces and $T_2$ uniform spaces
are denoted by ${\bold{T_2}}$, ${\bold{T_2Prox}}$ and ${\bold{T_2Unif}}$, 
respectively.
The categories of topological spaces, proximity spaces and uniform spaces
(without the $T_0,T_2,T_2$ axiom) 
are denoted by ${\bold{Top}}$, ${\bold{Prox}}$ and 
${\bold{Unif}}$, respectively. For any of these six categories (from 
${\bold{T_2}}$ till ${\bold{Unif}}$), $U$ will denote
their underlying set functors --- for which of these categories, that
will be given in our respective theorems. Also, if 
we have a subcategory ${\Cal{C}}$ of
them, $U$ may denote also the underlying set functor of ${\Cal{C}}$. It will
be always clear, which one do we mean. The category of $T_3$ topological
spaces is denoted by ${\bold{T_3}}$, and the notation $U$ will be used for it
in the above sense.

{\it{Subcategories are considered to be full and isomorphism closed, and
will be identified with the classes of their objects. A subcategory of}} 
${\bold{Top}}$, ${\bold{Prox}}$ {\it{or}} 
${\bold{Unif}}$ {\it{is non-trivial}} if it
contains a space with at least two points.

All the above six categories (from ${\bold{T_2}}$ till ${\bold{Unif}}$)
are complete, cocomplete, well-powered
and co-well-powered. Thus, epireflective subcategories can be characterized in
them, as those closed under products and extremal
subobjects. Also, epimorphisms can be factorized as a composition of a
bimorphism and an extremal epimorphism.
In ${\bold{T_2}}$, ${\bold{T_2Prox}}$, ${\bold{T_2Unif}}$ 
{\it{monomorphisms are the injections, 
epimorphisms are the dense maps, bimorphisms are the dense injections,
extremal monomorphisms are the closed embeddings, and 
extremal epimorphisms are the quotient maps in the respective categories}} 
(finest structures on surjective images making the surjective map a morphism).
In ${\bold{Top}}$, ${\bold{Prox}}$, ${\bold{Unif}}$
{\it{monomorphisms are the injections, 
epimorphisms are the surjections, bimorphisms are the bijections,
extremal monomorphisms are the embeddings, and 
extremal epimorphisms are the quotient maps in the respective categories}}.

For ${\bold{M}}$ a class of monomorphisms, or ${\bold{E}}$ a class of 
epimorphisms of some
category, we say that a subcategory ${\Cal{C}}$ is {\it{closed under 
${\bold{M}}$-subobjects}}, or is {\it{closed under ${\bold{E}}$-images}} if $C
\in $ ${\text{Ob}}\,{\Cal{C}}$ and $\exists m \in {\bold{M}}$, $m:D \to C$, or 
$\exists e \in {\bold{E}}$,  $e:
C \to D$ imply $D \in $ ${\text{Ob}}\,{\Cal{C}}$, respectively.

The {\it{covering character}} of a uniform space $X$, written as 
cov\,char\,$X$, is the smallest infinite
cardinal $\lambda $ such that $X$ admits a base
of uniform coverings of cardinalities less than $\lambda $. Equivalently, it
is the smallest infinite cardinal $\lambda $ such that $X$ has no uniformly
discrete subspace of cardinality $\lambda $ (\cite{I}). 
The completion of a uniform
space $X$ is denoted by $\gamma X$, and the precompact reflection of $X$
is denoted
by $pX$. Proximity spaces can be identified with precompact uniform spaces.
For general information about uniform spaces, cf. the book of Isbell \cite{I}.

\newpage

A concrete category 
$\langle {\Cal{C}}, U:{\Cal{C}} \to {\bold{Set}} \rangle $
with underlying set functor $U$
is called {\it{algebraic}} if ${\Cal{C}}$ has coequalizers,
$U$ has a left adjoint $F$, and $U$ preserves and reflects regular
epimorphisms. An algebraic category is {\it{varietal}} if additionally $U$
reflects congruence relations. An algebraic category is the same as a
{\it{quasivariety}}, i.e., all universal algebras of some given type
(infinitary operations allowed, which may form a proper class), 
closed under products and subalgebras. A
varietal category is the same as a {\it{variety}}, i.e., a quasivariety 
that is additionally closed under homomorphic images. Cf. \cite{HS07},
\cite{AHS}, or for a short description \cite{R82}.
We note that
for concrete categories $\langle {\Cal{C}}, U:{\Cal{C}} \to {\bold{Set}}
\rangle $
{\it{varietal}} is the same as
{\it{monadic}} or {\it{tripleable}} (cf. \cite{ML} and \cite{AHS}).

For category theory, we refer to \cite{ML}, \cite{H68}, \cite{HS07} and 
\cite{AHS}.


\head \S 2. Introduction\endhead


We begin with citing some theorems. The first one is a Birkhoff-type theorem
for ${\bold{T_2}}$.


\proclaim{Theorem A} {\rm{(D. Petz}} \cite{P}, {\rm{Theorem)}}.
Let ${\Cal{C}}$ be a subcategory of \,${\bold{T}}_2$. Then the following are
equivalent:

1) ${\Cal{C}}$ is an epireflective subcategory of \,${\bold{T}}_2$, 
closed under epi images. 

2) Either ${\text{\rm{Ob}}}\,{\Cal{C}}= 
\{ X \in {\text{\rm{Ob}}}\,{\bold{T_2}} \mid |X| \le 1 \} $, 
or \,${\text{\rm{Ob}}}\,{\Cal{C}}= 
\{ X \in {\text{\rm{Ob}}}\,{\bold{T_2}} \mid X$ is 
compact $T_2 \} $, or \,${\Cal{C}}= 
{\bold{T_2}}$.

If we write 1) as being closed under products, extremal
subobjects, bimorphic images and extremal epi images, then none
of these properties can be omitted, without invalidating the implication $1)
\Longrightarrow 2)$ of the theorem.
\endproclaim


We remark that \cite{P} did 
not decompose the hypotheses in the same way, therefore
we have to give examples for the last two properties. The 1-st, 2-nd, 4-th
properties are satisfied 
for $T_2$ spaces in which each at most countably infinite 
subset has a compact closure.
The 1-st, 2-nd, 3-rd properties are satisfied 
for $0$-dimensional compact $T_2$ spaces. 
(Both examples were given by \cite{P}).

One of the main results of G. Richter \cite{R89} 
was a generalization of Theorem A of Petz, with weaker (although more
complicated) hypotheses, which we give here.


\proclaim{Theorem B} {\rm{(Richter}} \cite{R89}, {\rm{Corollary 3.3)}}
Let ${\text{\rm{Ob}}}\,{\Cal{C}} \subset {\text{\rm{Ob}}}\,{\bold{T_2}}$. Then
the implication $1) \Longrightarrow 2)$
of Theorem A remains true under the following weaker
hypotheses. 
\newline
(P1'): the underlying set functor $U$ of ${\Cal{C}}$
has a left adjoint $F$ with pointwise dense unit $\eta : {\text{id}} _
{{\bold{Set}}} \to UF$ (i.e., for all 
$X \in {\text{\rm{Ob}}}\,{\bold{Set}}$ $\eta _X(X)$ is dense in $FX$);
\newline
(P2): ${\Cal{C}}$ is closed under surjective images;
\newline
(P3'): if $\eta _X: X \to UFX$ is a $C^*$-embedding (i.e., it underlies a 
$C^*$-embedding $X_d \to FX$, where $X_d$ is the discrete topological space
on the set $X$)
and $i:FX \to C$ is an
open, dense $C^*$-extension, such that in $C$ disjoint zero-sets can be
separated by a clopen set, then $C$ belongs to the maximal subcategory 
$\hat{{\Cal{C}}}$ of ${\bold{T}}_2$, containing ${\Cal{C}}$,
for whose underlying set functor $\hat{U}$ we have that
$F$ and $\eta $ still serve as left adjoint and unit of adjunction.
\endproclaim


Here (P1') is strictly weaker than epireflectivity, and (P3') is strictly
weaker than 

\newpage

closedness under dense extensions.
For more details we refer to \cite{R89} (in particular to its Theorem 3.1).

\vskip.1cm

H. Herrlich-G. R. Strecker \cite{HS71} 
initiated the investigation of algebraically
behaving subcategories of ${\bold{T}}_2$. 


\proclaim{Theorem C} {\rm{(H. Herrlich-G. R. Strecker}} \cite{HS71},
{\rm{Theorem, G. Richter}} \cite{R82},
{\rm{Corollary 1.6,}} \cite{R85a}, {\rm{Corollary 3.4)}}
Let ${\Cal{C}}$ be a subcategory of \,${\bold{T}}_2$.
Then the following are
equivalent:

1) ${\Cal{C}}$ is an epireflective subcategory of \,${\bold{T_2}}$, 
that is varietal. 

2)
Either ${\text{\rm{Ob}}}\,{\Cal{C}}= 
\{ X \in {\text{\rm{Ob}}}\,{\bold{T_2}} \mid |X| \le 1 \} $, 
or {\text{\rm{Ob}}}\,${\Cal{C}}= 
\{ X \in {\text{\rm{Ob}}}\,{\bold{T_2}} \mid X$ is 
compact $T_2 \} $.

If we write 1) as being closed under products, extremal
subobjects, and being varietal, then none
of these properties can be omitted, without invalidating the implication $1)
\Longrightarrow 2)$ of the theorem.
\endproclaim


For algebraic subcategories, G. Richter
proved an analogous result.


\proclaim{Theorem D} {\rm{(G. Richter}} \cite{R82}, {\rm{Corollaries 
1.5, 1.6)}}
Let \,${\Cal{C}}$ be a subcategory of \,${\bold{T_2}}$. 
Then the following are
equivalent:

1) ${\Cal{C}}$ is an epireflective subcategory of \,${\bold{T_2}}$ and is
algebraic.

2) ${\Cal{C}}$ is an epireflective subcategory 
of \,$\{ X \in {\text{\rm{Ob}}}\,{\bold{T_2}} \mid X$ 
is compact $T_2 \} $.
\endproclaim


If we write 1) as being closed under products, extremal
subobjects, and being algebraic, then none
of these properties can be omitted, without invalidating the implication $1)
\Longrightarrow 2)$ of the theorem. 
All, but the 1-st, 2-nd or 3-rd
property are
satisfied by examples 1) and 2) of \cite{HS71} and by 3) ${\bold{T_2}}$
(observe that its underlying set functor
does not reflect isomorphisms therefore it is not algebraic).
Example 1) is discrete topological spaces, 
which form even a varietal category \cite{HS71}.
Example 2) is the powers of a compact $T_2$
topological space, consisting of more than one point, that is strongly
rigid --- i.e., whose only continuous self-maps are the identity and the
constant maps (such spaces exist, cf. \cite{HS71}) --- 
together with the empty space, which category is even varietal,
cf. \cite{R92a}, p. 368.


Theorem C 
raises the analogous question for ${\bold{Top}}$. Both Theorem C, and its
word for word analogue for ${\bold{Top}}$, rather than ${\bold{T_2}}$, 
follow from the following theorem of Richter \cite{R91a}.


\proclaim{Theorem E} {\rm{(G. Richter}} \cite{R91a}, {\rm{Corollary 4.4)}}.
Let ${\Cal{C}}$ be a subcategory of \,${\bold{Top}}$. Then the following are
equivalent:

1) ${\Cal{C}}$
is closed under products, closed subspaces and is
varietal. 

2) Either ${\text{\rm{Ob}}}\,{\Cal{C}}= 
\{ X \in {\text{\rm{Ob}}}\,{\bold{Top}} \mid |X| \le 1 \} $, 
or ${\text{\rm{Ob}}}\,{\Cal{C}}= 
\{ X \in {\text{\rm{Ob}}}\,{\bold{Top}} \mid X$ is 
indiscrete$ \} $,
or ${\text{\rm{Ob}}}\,{\Cal{C}}= 
\{ X \in {\text{\rm{Ob}}}\,{\bold{Top}} \mid X$ is 
compact $T_2 \} $.
\endproclaim


Here none of the properties in 1) can be omitted,
without invalidating the implication $1) \Longrightarrow 2)$ of the theorem,
as follows from Theorem C and from the fact that $\bold{T_2}$ is closed under
products and closed subspaces.

Earlier, Richter \cite{R85a}, 
Corollary 3.4 proved a weaker result. Namely, he added
to the hypotheses of Theorem E
that the two-point discrete space belongs to ${\text{\rm{Ob}}}\,{\Cal{C}}$, 
cf. \cite{R85a}, p. 80
(then, of course, in 2) only the third case is possible.) 


\vskip.1cm

One could obtain a 
common generalization (roughly) of Theorems A and 1 (cf. \S 3, dealing with
epireflective subcategories of ${\bold{Top}}$, closed under bijective images),

\newpage

in the following way.
We consider 
subcategories of ${\bold{Top}}$ that are productive, are closed only 
under closed subspaces (like in Theorem A), and only under surjective images
(like in a weakened variant of Theorem 1, replacing ``bimorphic'' with
``epimorphic''). 
Such a theorem is available, but only for $T_3$ spaces.


\proclaim{Theorem F} {\rm{(Kannan-Soundararajan}} \cite{KS}, {\rm{Theorem)}}
Let ${\Cal{C}}$ be a subcategory 
of ${\bold{T_3}}$.
Then the following are equivalent. 

1) ${\Cal{C}}$ is productive,
closed-hereditary, and is closed under surjective images. 

2) Either 

A) {\rm{Ob}}\,${\Cal{C}}=\{ X \mid X$ is \,$T_3$ and \,\,$|X| \le 1 \} $, or

B) there exists a
class ${\Cal{F}}$ 
of ultrafilters $p$ on some sets $S_p$, such that {\rm{Ob}}\,${\Cal{C}}$ 
consists of all \,$T_3$ spaces $X$ satisfying the following property. For
$p \in {\Cal{F}}$ and $f:S_p \to X $ any function
there exists a continuous
extension of $f$, namely ${\overline{f}}:S_p \cup \{ p \} \to X$, 
where $S_p \cup \{ p \} $ has
the subspace topology inherited from $\beta [(S_p)_d]$, and
where $(S_p)_d$ is the discrete space on $S_p$.
\endproclaim


Here none of the properties in 1) can be omitted,
without invalidating the implication $1) \Longrightarrow 2)$ of the theorem.
This is shown by the following
examples in ${\bold{T_3}}$: 
finite spaces, connected spaces, zero-dimensional compact $T_2$
spaces (the 1-st and 3-rd examples are taken from \cite{P}).

\definition {Remark} 
Let us exclude the trivial case ${\text{Ob}}\,{\Cal{C}} = 
\{ X \mid X {\text{ is }} T_3 {\text{ and }}
|X| \le 1 \} $. Then the underlying set functor $U:{\Cal{C}}
\to {\bold{Set}}$ has a left adjoint $F$ that has a natural transformation to  
the functor
$X \to \beta (X _d)$, with all components embeddings (\cite{KS}, essentially
Step 4, p. 143, applied to $T_3$ spaces). 
Then we can recover a
(maximal) class ${\Cal{F}}$, by considering the spaces $FX$, for all $X \in
{\text{Ob}}\,{\bold{Set}}$: the points of all these spaces
$FX$ will give the class of ultrafilters mentioned in Theorem F.
(This is the construction of \cite{KS}, Step 6,
p. 144, except that there fixed ultrafilters are not considered, but that does
not change matters). 
Thus we have for this (maximal) class ${\Cal{F}}$ that for
$X \in {\text{Ob}}\,{\bold{Set}}$, and for any ${\bold{Set}}$-morphism
$f: Y \to X$ there exists a domain-codomain  
extension of $f$ to a ${\bold{T_3}}$-morphism 
$Ff:FY \to FX$, i.e., $(Ff)(FY) \subset FX$.
In short: ``${\Cal{F}}$ is closed under images''.
(This is a special case of the statement of Theorem F, 2), B)
applied to the space
$FX$ rather than $X$ in Theorem F 2) B)). 
Observe that \cite{KS} Theorem did not contain this
property of ${\Cal{F}}$ explicitly. 
In fact this property is necessary (and sufficient) for $F$ to be a
left adjoint of $U$ even when restricted to the minimal class 
$\{ FX \mid X \in {\text{Ob}}\,{\bold{Set}} \} $ 
(Kleisli adjunction, \cite{AHS}, 20.39, 20.B).

On the other hand, for given ${\Cal{F}}$, the class ${\Cal{C}}$ constructed in
Theorem F is a maximal subclass of ${\bold{T_3}}$ for which $F$ and $\eta $ are
left adjoint to $U$ and unit of adjunction.

Recall that all subcategories considered are isomorphism closed. Observe also
that $F$ and $\eta $ are defined only up to isomorphisms, forming respective
commutative diagrams. We can eliminate these ambiguities by using some specific
construction of $\beta (X_d)$, e.g., with ultrafilters, and considering 
$$
\eta_X X \subset UFX {\text{ and }} FX \subset \beta (X_d)
\tag *
$$
(thus the usual embeddings are realized by embeddings of subsets/subspaces).

Let us denote, for given (maximal) ${\Cal{C}}$ and (maximal) 
${\Cal{F}}$, by ${\Cal{C}}({\Cal{F}})$ and
${\Cal{F}}({\Cal{C}})$ the (maximal) subcategory ${\Cal{C}}$
constructed for ${\Cal{F}}$ in Theorem F, 2), B) and the (maximal) class
${\Cal{F}}$ constructed in the proof of \cite{KS}, Theorem, Step 6.

\newpage

\cite{KS} did not completely clarify the situation. Namely, for a category
${\Cal{C}}$ there exists a class ${\Cal{F}}$ of ultrafilters making 2) B) of 
Theorem F true. However the questions, 
which classes ${\Cal{F}}$ of ultrafilters arise this way, and
possibly when are the corresponding categories ${\Cal{C}}$
equal, are not considered there.
As already mentioned, the class ${\Cal{F}}$ is ``closed under images'', so
we need to consider this question only for classes of ultrafilters 
``closed under images''. 

The class of the spaces $X$ described in Theorem F, 2), B) is the largest class
${\Cal{C}}_{\text{max}}({\Cal{F}})$ (in ${\bold{T_3}}$!) 
for which $F$ and $\eta $ are left adjoint to $U$ and unit of
adjunction. This ${\Cal{C}}_{\text{max}}({\Cal{F}})$ determines uniquely 
$F$ and $\eta $ as left adjoint and unit of adjunction, by convention
\thetag{*}. 

Similarly, the class ${\Cal{F}}_{\text{max}}({\Cal{C}})$ gives exactly 
the Kleisli
adjunction (minimal adjunction) associated to the adjunction $F \dashv U$.

Beginning with the Kleisli adjunction, then taking the maximal (in
${\bold{T_3}}$) 
adjunction with given $F$ and $\eta $, and turning once more to the
Kleisli adjunction clearly gives back the original Kleisli adjunction (by
\thetag{*}).

Beginning with the maximal (in ${\bold{T_3}}$) adjunction with given $F$ and
$\eta $, then taking the Kleisli adjunction, and turning once more to the
maximal (in ${\bold{T_3}}$) adjunction with given $F$ and $\eta $ clearly gives
back the original maximal adjunction.

This settles the case of maximal subcategories
${\Cal{C}}_{\text{max}}({\Cal{F}})$. It will suffice to prove that under 1)
of Theorem F each subcategory ${\Cal{C}}$ is maximal.

Observe that the proof of \cite{KS}, Theorem, Step 7 in fact
proves the following. Let $C \in {\text{Ob}}\,{\Cal{C}}_{\text{max}}
({\Cal{F}})$. Then
$C \in {\text{Ob}}\,{\Cal{C}}$. (Namely there $C$ is the surjective image of 
$FUC$, by $\varepsilon _C$, the counit of the adjunction.) This proves
${\Cal{C}}_{\text{max}}({\Cal{F}}) \subset {\Cal{C}}$, i.e., that each  
subcategory ${\Cal{C}}$ in 1) of Theorem F is maximal.

That is, we have shown the following addition to Theorem F.
\enddefinition

\proclaim{Proposition}
In Theorem F, 2), B), we may additionally suppose that ${\Cal{F}}$ is "closed
under images" {\rm{(definition cf. above)}}. Under this restriction, Theorem
F, 2), B) 
establishes a bijection between the subcategories satisfying Theorem F,
1), and the classes ${\Cal{F}}$ of ultrafilters "closed under images".
$\blacksquare $
\endproclaim

If in Theorem F 1) 
we write instead of closedness under surjective images closedness
under bijections and closedness under
extremal epi images, then $0$-dimensional compact $T_2$
spaces are epireflective and
closed under bijective images in ${\bold{T_3}}$, but are not of the form in 2). 

\definition{Problem 1}
Find a fourth
example (if it exists)
that is epireflective and closed under extremal epi images, 
but is not of the form in 2) (the proof in \cite{KS}
seems to use, by $\varepsilon _C$, 
both closedness under bijections and extremal epi images).
\enddefinition

Extensions of Theorem F cf. in the paper of Hager \cite{Ha}, to the case of
a concrete category. say, over ${\bold{Set}}$. Then his theorem is specialized
to ${\bold{T_2Prox}}$ and ${\bold{T_2Unif}}$. (And also for $T_2$ cozero
spaces, where a {\it{$T_2$ cozero space}} 
can be easiest defined as the cozero sets
of all uniformly continuous real valued 
functions for some $T_2$ uniformity on the underlying set, and a cozero
morphism is a set morphism, such that the inverse image of a cozero set is
also a cozero set. More about this cf. in \cite{Ha}.)
For 

\newpage

details we have to refer to \cite{Ha}.

\vskip.1cm

Richter \cite{R80/81}, \cite{R82}, \cite{R85a}, \cite{R85b}, 
\cite{R89}, \cite{R91a}, \cite{R91b}, \cite{R92a}, \cite{R92b}, \cite{R99} 
contain much related material.

\vskip.1cm

We will prove analogues of these theorems for ${\bold{Top}}$,
${\bold{Prox}}$, 
${\bold{T_2Unif}}$, ${\bold{Unif}}$. 


\head \S 3. Theorems\endhead

The first three theorems will deal with epireflective categories closed under
epimorphic or bimorphic images.


\vskip.1cm

First we give a simple proof of an analogue of Theorem A for ${\bold{Top}}$,
with less hypotheses. 


\proclaim{Theorem 1}
Let ${\Cal{C}}$ be a subcategory of \,${\bold{Top}}$. Then the following are
equivalent:

1) ${\Cal{C}}$ an epireflective subcategory of \,${\bold{Top}}$, 
closed under bimorphic images. 

2) 
Either ${\text{\rm{Ob}}}\,{\Cal{C}}= 
\{ X \in {\text{\rm{Ob}}}\,{\bold{Top}} \mid |X| \le 1 \}$, 
or ${\text{\rm{Ob}}}\,{\Cal{C}}= 
\{ X \in {\text{\rm{Ob}}}\,{\bold{Top}} \mid X$ is 
indiscrete$ \} $, or \,${\Cal{C}}= {\bold{Top}}$.

If we write 1) as being closed under products, extremal
subobjects and bimorphic images, then none
of these properties can be omitted, without invalidating 
the implication $1) \Longrightarrow 2)$ of the theorem.
\endproclaim


Next we give the analogue of Theorem A for ${\bold{T_2Unif}}$.


\proclaim{Theorem 2}
Let ${\Cal{C}}$ be a subcategory of \,${\bold{T_2Unif}}$. 
Then the following are equivalent:

1) ${\Cal{C}}$ is an epireflective subcategory of \,${\bold{T_2Unif}}$, closed
under epi images. 

2)
Either ${\text{\rm{Ob}}}\,{\Cal{C}}= 
\{ X \in {\text{\rm{Ob}}}\,{\bold{T_2Unif}} \mid |X| \le 1 \} $, 
\,or ${\text{\rm{Ob}}}\,{\Cal{C}}= 
\{ X \in {\text{\rm{Ob}}}\,{\bold{T_2Unif}} \mid X$ is 
compact $T_2 \} $, 
or there exists an infinite cardinal $\lambda _0$, such that
${\text{\rm{Ob}}}\,{\Cal{C}} = 
\{ X \in {\text{\rm{Ob}}}\,{\bold{T_2Unif}} \mid X$ has a covering
character at most $\lambda _0 \}$, or \,${\Cal{C}}={\bold{T_2Unif}}$.

If we write 1) as being closed under products, extremal
subobjects, bimorphic images and extremal epi images, 
then none
of these properties can be omitted, without invalidating the implication $1)
\Longrightarrow 2)$ of the theorem.
\endproclaim


Next we turn to a common analogue of Theorems 1 and 2, for ${\bold{Unif}}$.


\proclaim{Theorem 3}
Let ${\Cal{C}}$ be a subcategory of \,${\bold{Unif}}$. Then the following are
equivalent:

1) ${\Cal{C}}$ is an epireflective subcategory of \,${\bold{Unif}}$, closed
under bimorphic images.

2) Either ${\text{\rm{Ob}}}\,{\Cal{C}}= 
\{ X \in {\text{\rm{Ob}}}\,{\bold{Unif}} \mid |X| \le 1 \} $, 
or ${\text{\rm{Ob}}}\,{\Cal{C}}= 
\{ X \in {\text{\rm{Ob}}}\,{\bold{Unif}} \mid X$ is indiscrete$\} $, 
or there exists an infinite cardinal $\lambda _0$
such that \,${\text{\rm{Ob}}}\,{\Cal{C}}= 
\{ X \in {\text{\rm{Ob}}}\,{\bold{Unif}} \mid X$ has a covering
character at most $\lambda _0 \}$, or \,${\Cal{C}}={\bold{Unif}}$.

If we write 1) as being closed under products, extremal
subobjects and bimorphic images, then none
of these properties can be omitted, without invalidating the implication $1)
\Longrightarrow 2)$ of the theorem.
\endproclaim


\definition{Problem 2}
There arises the question about the uniform version of Theorem F of
Kannan-Soundararajan. That is, we suppose closedness under
products, closed subspaces and surjective images. This would be a common
generalization of Theorems 2 and 3.
As mentioned after Theorem F, Theorems 2 and 3 have a common generalization
in Hager \cite{Ha}. However, the
description in \cite{Ha} does not seem to imply in an evident way our concrete
descriptions in our Theorems 2 and 3. (It uses for the description
some class of epimorphisms 
of the category, of which there are illegitimely many; and it does not seem to
be evident how to identify these classes and concretize the description
in our concrete cases.)
Also, our proofs are independent of \cite{Ha}. 

\newpage

However this seems to
be a question much more complicated than Theorem F of Kannan-Soundararajan.
Let us restrict our attention to the case of ${\bold{T_2Unif}}$. 
Of course, we have as examples all $T_2$ uniform spaces, whose
underlying topological spaces form a class (of Tychonoff spaces!)
as described in Theorem
F. Moreover, if we take a cardinal $\lambda _0 \ge \aleph _0$, we have an
example all those uniform spaces, whose
underlying topological spaces form a class as described in Theorem F, and
whose covering characters are at most $\lambda _0$. (For 
$\lambda _0 = \aleph _0$ we have examples for proximities.) The problem is
again that we cannot identify the class of epimorphisms whose existence is
stated in \cite{Ha} and cannot concretize the description in 
\cite{Ha}. So a concrete, usable description still is missing.
\enddefinition


The next three theorems will deal with epireflective subcategories, or
subcategories closed under products and closed subspaces,
which are algebraic or varietal. 


\proclaim{Theorem 4}
Let ${\bold{T}}=
{\bold{Top}}$ or \,${\bold{T}}={\bold{Prox}}$ or \,${\bold{T}}={\bold{Unif}}$.
Let ${\Cal{C}}$ be a subcategory of \,${\bold{T}}$. Then the following are
equivalent:

1) 
${\Cal{C}}$ is an epireflective subcategory of \,${\bold{T}}$, 
that is algebraic. 

2)
Either ${\text{\rm{Ob}}}\,{\Cal{C}}= 
\{ X \in {\text{\rm{Ob}}}\,{\bold{T}} \mid |X| \le 1 \} $, 
or ${\text{\rm{Ob}}}\,{\Cal{C}}
= \{ X \in {\text{\rm{Ob}}}\,{\bold{T}} \mid X$ is indiscrete$\} $.

If we write 1) as being closed under products, extremal
subobjects, and being algebraic, then none
of these properties can be omitted, 
without invalidating the implication $1) \Longrightarrow 2)$
of the theorem.
\endproclaim

Although the cases ${\bold{T}}={\bold{Prox}}$ and ${\bold{T}}={\bold{Unif}}$
of Theorem 4 are covered by the next
theorem (namely
the minimal non-trivial epireflective subcategory of compact $T_2$
proximity or uniform spaces is that of $0$-dimensional compact $T_2$ 
proximity or uniform
spaces, and its
hereditary hull, taken 
in ${\bold{Prox}}$ or ${\bold{Unif}}$, 
contains also non-compact proximity or uniform spaces),
its proof in Theorem 4 is much simpler than the proof of 
Theorem 5.


\proclaim{Theorem 5}
Let ${\Cal{C}}$ be a subcategory of \,${\bold{Unif}}$.
Then the following are equivalent: 

1) ${\Cal{C}}$ is closed under products and closed subspaces and is algebraic. 

2) Either
${\text{\rm{Ob}}}\,{\Cal{C}}= \{ X \in {\text{\rm{Ob}}}\,{\bold{Unif}} \mid X$ 
is indiscrete$\} $,
or \,${\Cal{C}}$ is an epireflective subcategory of 
$ \{ X \in {\text{\rm{Ob}}}\,{\bold{Unif}} \mid X$ is compact $T_2 \} $. 

None
of the properties of 1) of this theorem
can be omitted, without invalidating the implication
$1) \Longrightarrow 2)$ of the theorem. 
\endproclaim


Next we turn to an analogue of Theorem E for ${\bold{Unif}}$. This theorem
implies the description of epireflective and varietal subcategories both in
${\bold{T_2Unif}}$ and ${\bold{Unif}}$.


\proclaim{Theorem 6} 
Let ${\Cal{C}}$ be a subcategory of \,${\bold{Unif}}$. Then the following are
equivalent:

1) ${\Cal{C}}$ is a subcategory of \,${\bold{Unif}}$, 
closed under products and closed subspaces, that
is varietal.

2)
Either ${\text{\rm{Ob}}}\,{\Cal{C}}= 
\{ X \in {\text{\rm{Ob}}}\,{\bold{Unif}} \mid |X| \le 1 \} $, 
or ${\text{\rm{Ob}}}\,{\Cal{C}}= 
\{ X \in {\text{\rm{Ob}}}\,{\bold{Unif}} \mid X$ is 
indiscrete$ \} $, 
or ${\text{\rm{Ob}}}\,{\Cal{C}}= 
\{ X \in {\text{\rm{Ob}}}\,{\bold{Unif}} \mid X$ is 
compact $T_2 \} $. 

None of the properties in 1) can be omitted, 
without invalidating the implication $1) \Longrightarrow 2)$ of the theorem.
\endproclaim


\definition{Problem 3}
What remains open, is the following question, that would include Theorems D,
E and 4 (for ${\bold{Top}}$), 
and would be an analogue of Theorem 5. Namely, can one describe all
subcategories ${\Cal{C}}$ 
of ${\bold{Top}}$, closed under products and closed 

\newpage

subspaces, 
and
being algebraic? Are there more such subcategories than described in Theorems
D and E? 
(Observe that the uniform case is settled by Theorem 5. 
Thus the situation
is just the converse of the situation mentioned in Problem 2,
where the uniform case seems to be much more complicated.)
\enddefinition


\head \S 4. Proofs\endhead


\demo{Proof of Theorem 1}
We only need to prove $1) \Longrightarrow 2)$.

The empty product, i.e., the one-point space belongs to 
${\text{\rm{Ob}}}\,{\Cal{C}}$, as
well as its (closed) subspace the empty set. Hence, 
$\{ X \in {\text{\rm{Ob}}}\,{\bold{Top}} \mid |X| \le 1 \} \subset
{\text{\rm{Ob}}}\,{\Cal{C}}$. If here we have equality, we are done.

Therefore we may suppose that some space $X$ belongs to 
${\text{\rm{Ob}}}\,{\Cal{C}}$,
where $|X| \ge 2$. Then any of its two-point subspaces belongs to 
${\text{\rm{Ob}}}\,{\Cal{C}}$ as well, 
hence we may suppose $|X|=2$. Then its bijective
image the two-point indiscrete subspace belongs to 
${\text{\rm{Ob}}}\,{\Cal{C}}$ as well,
hence we may suppose that $X$ is the two-point indiscrete space $I_2$.
Then any subspace of any power of $I_2$ belongs to 
${\text{\rm{Ob}}}\,{\Cal{C}}$, hence
$\{ X \in {\text{\rm{Ob}}}\,{\bold{Top}} \mid X$ is indiscrete$ \} \subset
{\text{\rm{Ob}}}\,{\Cal{C}}$. If here we have equality, we are done.

Therefore we may suppose that some space $X$ belongs to 
${\text{\rm{Ob}}}\,{\Cal{C}}$,
where $X$ is non-indiscrete. Then $X$ has a non-indiscrete two-point subspace,
that of course belongs to ${\text{\rm{Ob}}}\,{\Cal{C}}$, 
hence we may suppose $|X|=2$.
Then the Sierpi\'nski space is a bijective image of $X$, therefore it belongs
to ${\text{\rm{Ob}}}\,{\Cal{C}}$ as well. 
Since any $T_0$ topological space is a subspace of
a power of the Sierpi\'nski space, hence
$\{ X \in {\text{\rm{Ob}}}\,{\bold{Top}} \mid X$ is $T_0 \} \subset
{\text{\rm{Ob}}}\,{\Cal{C}}$. Finally, any topological space is a subspace 
of a product of a
$T_0$ space and an indiscrete space. Hence ${\Cal{C}} = {\bold{Top}}$.

There remains to give three examples. 
All but the 1-st, 2-nd, or 3-rd properties are satisfied by
the subclasses of ${\bold{Top}}$ consisting
of finite spaces, of connected spaces (both being closed even
under all surjective images), or of $T_0$ spaces,
respectively.
$ \blacksquare $
\enddemo


We begin the proof of Theorem 2
with a simple lemma, that is known. For 1) of Lemma 1 (for realcompact
spaces), cf. \cite{GJ}, Theorem 8.9, and for 2) of Lemma 1
(also for realcompact spaces),
cf. \cite{GJ}, Theorem 8.13. A categorical generalization of both 1) and 2), 
namely that epireflective
subcategories are strongly closed under limits, with an explanation that 1)
and 2) are particular cases of this general statement, cf. in \cite{H68}, 
\S 9.3. We state our Lemma for ${\bold{T_2Unif}}$.


\proclaim{Lemma 1} {\rm{(\cite{GJ}, \cite{H68}, cited just before this Lemma)}}
Let ${\Cal{E}}$ be an epireflective subcategory of \,${\bold{T_2Unif}}$. 

1) Let $X \in {\text{\rm{Ob}}}\,{\bold{T_2Unif}}$, and let 
$X_{\alpha }$,
for $\alpha \in A$, be subspaces of $X$,
such that for each $\alpha \in A$ we have $X_{\alpha }
\in $ {\text{\rm{Ob}}}\,${\Cal{E}}$. Then $\cap _{\alpha \in A} X_{\alpha }
\in $ {\text{\rm{Ob}}}\,${\Cal{E}}$.

2) Let $X \in {\text{\rm{Ob}}}\,{\Cal{E}}$, 
$Y \in {\text{\rm{Ob}}}\,{\bold{T_2Unif}}$, $Z \subset Y$, $Z \in 
{\text{\rm{Ob}}}\,{\Cal{E}}$ and let $f:X \to Y$ be uniformly
continuous. Then $f^{-1}(Z) \in {\text{\rm{Ob}}}\,{\Cal{E}}$.
$ \blacksquare $
\endproclaim


Next we give a certain
uniform analogue of well-known theorems for topological spaces, cf.
[E], Exercises 4.2.D and 
4.4.J and Theorem 4.4.15, about representing topological, or metric spaces
as images of certain spaces under certain types of mappings. In 
particular, these statements characterize the class of
first countable $T_0$ spaces, or metric spaces, as the open, or 
perfect images of 
subspaces of Baire 

\newpage

spaces $D_{\lambda }^{\aleph
_0}$ --- where $D_{\lambda }$ is a discrete topological space of cardinality
$\lambda $, and where
$\lambda $ equals the weight of the space to be represented --- respectively.
Some more specialized theorems of this type cf., e.g., in \cite{M}.
(About inverse limits of uniform spaces, to be used in the proof of Lemma 2,
cf. \cite{I}, 
\S IV, subchapter ``Inverse limits''.) 


\proclaim{Lemma 2}
Let $M$ be a complete metric space with covering character at most $\lambda
_0$, where $\lambda _0$ is an infinite cardinal. Then there is a dense,
uniformly continuous map from a closed subspace of a countable product (taken
in ${\bold{T_2Unif}}$) of discrete uniform spaces, of cardinalities less than
$\lambda _0$, to $M$.
\endproclaim


\demo{Proof}
Let $(M, \varrho )$ be our complete metric space, with cov\,char\,$X \le \lambda
_0$. By replacing the original metric $\varrho $ 
by $(1 - \varepsilon ) \varrho / (1
+ \varrho )$, if necessary, we may assume that diam\,$M < 1$. For each
integer $n \ge 0$
we will define sets $M_n \subset M$ as follows. $M_n$ is a maximal
subset of $M$, containing $M_{n-1}$ (for $n=0$ 
we let $M_{-1} = \emptyset $),
such that any two different points of $M_n$ have a distance at least
$1/2^n$. Clearly $|M_0|=1$, and for all $n$ we have $|M_n| < \lambda
_0$. (This is true also for cov\,char\,$M = \aleph _0$, i.e., when $M$ is
precompact.)

For $n \ge 0$ we define maps $f_n : M_{m+1} \to M_n$, such that $f_n$ is
identity on $M_n$, and else, for $m_{n+1} \in M_{n+1} \setminus M_n$, we have 
$$
\varrho (m_{n+1}, f_n(m_{n+1})) < 1/2^n .
\tag *
$$ 
The existence of $f_n(m_{n+1})$ follows from the
maximality property of $M_n$. Of course, inequality \thetag{*} holds for all
$m_{n+1} \in M_{n+1}$. The same maximality property, for each $n$, implies that
$\cup _{n=0}^{\infty } M_n$ is dense in $M$.  
  
We define a partial order $\le $ on $\cup _{n=0}^{\infty } M_n$, as the 
transitive (and reflexive) hull of the relation 
$$
\{ ( f_n(m_{n+1}), m_{n+1})
\mid n \ge 0,\,\,m_{n+1} \in M_{n+1} \} .
$$ 
This gives a tree structure on
$\cup _{n=0}^{\infty } M_n$, and any two points of 
$\cup _{n=0}^{\infty } M_n$ have a greatest lower bound. The $0$-th, $1$-st,
$2$-nd, $\dots $ levels of the tree are $M_0, M_1\setminus M_0, M_2 \setminus
M_1, \dots $\,.

Then $M_0 \overset{f_0}\to \longleftarrow 
M_1 \overset{f_1}\to \longleftarrow \dots $ forms an inverse system of
complete, and in fact, uniformly discrete uniform spaces. Its inverse limit 
$\underleftarrow{\lim} (M_n,f_n)$ is a complete uniform space, and $\cup
_{n=0}^{\infty } M_n$ has a natural embedding $i$ to  
$\underleftarrow{\lim} (M_n,f_n)$:
to $m_n \in M_n$ we let correspond the thread (branch)
$$
\cases
i(m_n) := \langle f_0f_1\dots f_{n-1}(m_n), f_1\dots f_{n-1}(m_n), \dots , \\
f_{n-2}f_{n-1}
(m_n),f_{n-1}(m_n),m_n,m_n, m_n, \dots \rangle .
\endcases
$$ 
We define a metric
$d$ on $i(\cup_{n=0}^{\infty } M_n)$ as follows. For $m_{n_1} \in M_{n_1}
\setminus M_{n_1-1} $ and $m_{n_2} \in M_{n_2}
\setminus M_{n_2-1}$, we let $d(im_{n_1},im_{n_2}) := 1/2^n$,
where the greatest lower bound of $m_{n_1}$ and $m_{n_2}$ is on the $n$-th
level, where $n \ge 0$. 

This can be extended to a metric $d$ on 
$\underleftarrow{\lim} (M_n,f_n)$ as follows. The distance of two different
threads (branches) is $1/2^n$, if the threads are identical exactly on the 
$0$-th, $1$-st, $2$-nd, $\dots , n$-th levels. (This metric is
non-Archimedean, i.e., we have $d(x,z) \le \max \{ d(x,y), d(y,z) \} $, thus, 
in particular,
$d(im_1,im_3) \le \max \{ d(im_1,im_2), 
d(im_2,$

\newpage

$im_3) \} $.) Then $i(\cup _{n=0}^{\infty } M_n)$ is dense in 
$\underleftarrow{\lim} (M_n,f_n)$. Observe that $\underleftarrow{\lim}
(M_n,f_n)$ is a closed
subspace of the product $\prod _{n=0}^{\infty } M_n$.

Let us map $i(\cup_{m=0}^{\infty }M_n)$ to $\cup_{m=0}^{\infty }M_n$ 
by the left
inverse $j$ of the embedding $i$, when $i$ is considered here as a map from 
$\cup_{m=0}^{\infty }M_n$ to $i(\cup_{m=0}^{\infty }M_n)$. 
We assert that $j$ is a
Lipschitz map with Lipschitz constant $4$. In fact, let $m_{n_1} \in 
M_{n_1} \setminus M_{n_1-1}$ and $m_{n_2} \in 
M_{n_2} \setminus M_{n_2-1}$ have greatest lower bound on level $n$; thus 
$d(im_{n_1}, im_{n_2})=1/2^n$. Then 
$$
\cases
\varrho (m_{n_1}, m_{n}) \le 
\varrho \left( m_{n_1},f_{n_1-1}(m_{n_1}) \right) + 
\varrho \left( f_{n_1-1}(m_{n_1}),f_{n_1-2}f_{n_1-1}(m_{n_1}) \right) + 
\\
\dots +\varrho \left( f_{n+1}f_{n+2} \dots f_{n_1-2}f_{n_1-1}(m_{n_1}),
f_nf_{n+1}f_{n+2} \dots f_{n_1-2}f_{n_1-1}(m_{n_1}) \right) 
\\
<1/2^{n_1-1}+1/2^{n_1-2}+ \dots +1/2^n < 2/2^n.
\endcases
$$ 
Similarly, 
$\varrho (m_{n_2}, m_{n}) < 2/2^n$, hence 
$$
\varrho (m_{n_1}, m_{n_2}) < 4/2^n = 4d(im_{n_1},im_{n_2}),
$$ 
as claimed above.

Now recall that 
$i(\cup _{n=0}^{\infty } $
$M_n)$ is dense in the complete metric space
$\underleftarrow{\lim} (M_n,f_n)$, 
and $\cup _{n=0}^{\infty } M_n$ 
is dense in the complete metric space
$M$. 
Then $j$ has an extension $\varphi :
\underleftarrow{\lim} (M_n,f_n) \to M$, that is Lipschitz with constant $4$,
hence is a uniformly continuous and dense map. 
$\blacksquare $
\enddemo


\demo{Proof of Theorem 2}
We only need to prove $1) \Longrightarrow 2)$.

{\bf{1.}}
Like in the proof of Theorem 1, second paragraph, we see that 
$\{ X \in {\text{{\rm{Ob}}}}\,{\bold{Unif}} \mid |X| \le 1 \} \subset $
${\text{\rm{Ob}}}\,{\Cal{C}}$. If here we have equality, we are done.

Now suppose that here we do not have equality, i.e., ${\Cal{C}} $ contains a
$T_2$ uniform space $X$ with at least two points. Then $X$ has a closed
subspace consisting of two points, i.e., the discrete two-point space $D_2$,
that therefore
belongs to ${\text{\rm{Ob}}}\,{\Cal {C}}$. 
Then all closed subspaces of all finite
powers of $D_2$ belong to ${\text{\rm{Ob}}}\,{\Cal {C}}$, hence 
$\{ X \in {\text{{\rm{Ob}}}}\,{\bold{T_2Unif}} \mid X $ is a finite discrete 
space$ \} \subset {\text{\rm{Ob}}}\,{\Cal{C}}$.

Also $D_2^{\aleph _0} \in {\text{\rm{Ob}}}\,{\Cal {C}}$ (power meant in
${\bold{T_2Unif}}$), i.e., the Cantor set with its
unique compatible uniformity belongs to
${\text{\rm{Ob}}}\,{\Cal{C}}$.
Then also its uniformly continuous image $[0,1]$ belongs to 
${\text{\rm{Ob}}}\,{\Cal{C}}$,
and all its powers $[0,1]^{\alpha }$ belong to 
${\text{\rm{Ob}}}\,{\Cal{C}}$, as well as all
their closed  subspaces. That is, 
$\{ X \in {\text{{\rm{Ob}}}}\,{\bold{T_2Unif}} \mid X$ 
is compact $T_2 \} \subset {\text{\rm{Ob}}}\,{\Cal{C}}$. 
If here we have equality, we are done.
(Up to this point, the proof is essentially the same, as in \cite{W}, p. 51,
\cite{HS71} and \cite{P}.)

Now suppose that here we do not have equality, 
i.e., ${\text{\rm{Ob}}}\,{\Cal{C}} $ contains a
$T_2$ uniform space $X$ that is not compact. Then its uniformly continuous
image $pX$, its precompact reflection, is homeomorphic to $X$, hence also is
non-compact, and belongs to 
${\text{\rm{Ob}}}\,{\Cal{C}}$. Thus we may assume that $X \in 
{\text{\rm{Ob}}}\,{\Cal{C}}$ is precompact, non-compact. 
Then $\gamma X$, its completion,
is a proper superset of $X$. Further, $\gamma X$ is compact $T_2$.

Following \cite{P}, choose $a \in X$ and $b \in \gamma X \setminus X$. 
Then $\{ a,b \} ^{\aleph _0} \subset
(\gamma X)^{\aleph _0}$, and $\{ a,b \} ^{\aleph _0}$ is the Cantor set with
its unique compatible uniformity, $C$, say.
By $X \in {\text{\rm{Ob}}}\,{\Cal{C}}$ 
we have $X^{\aleph _0} \in {\text{\rm{Ob}}}\,{\Cal{C}}$, hence
also any subspace of $(\gamma X)^{\aleph _0}$, containing $X^{\aleph _0}$
(that is dense in $(\gamma X)^{\aleph _0}$), 
belongs to ${\text{\rm{Ob}}}\,{\Cal{C}}$. 
In particular, $(\gamma X)^{\aleph _0} \setminus
\{ \langle b,b, \dots \rangle \} \in {\text{\rm{Ob}}}\,{\Cal{C}}$. 
This last subspace has
as closed subspace $\{ a,b \} ^{\aleph _0} \setminus 
\{ \langle b,b, \dots \rangle \} $. Hence, using for $C$ be the usual ternary
representation of the Cantor set, we have $C \setminus \{ 0 \} \in 
{\text{\rm{Ob}}}\,{\Cal{C}}$. Then, 
for the usually constructed surjection $f:C \to 
[0,1]$, we
have $f(C \setminus \{ 0 \} ) =$

\newpage

$(0,1] \in {\text{\rm{Ob}}}\,{\Cal{C}}$.

{\bf{2.}}
The class of cardinalities $\lambda $ (finite of infinite), 
for which the discrete space $D_{\lambda
}$ of cardinality $\lambda $ belongs to ${\text{\rm{Ob}}}\,{\Cal{C}}$, 
forms an initial
segment of all cardinalities, i.e., it is of the form $\{ \lambda \mid \lambda
< \lambda _0 \} $ or $\{ \lambda \mid \lambda $ is a cardinal$\} $. In the
first case, by the last sentence of the second paragraph of
{\bf{1}}, we have $\lambda _0 \ge \aleph _0$.

{\bf{3.}}
We begin with the case when this initial segment is $\{ \lambda \mid \lambda
< \lambda _0 \} $. No $T_2$ uniform space in 
${\text{\rm{Ob}}}\,{\Cal{C}}$ can have a
covering character greater than $\lambda _0$, since such a space contains a
closed subspace $D_{\lambda _0}$, and then we would have $D_{\lambda _0} \in
{\text{Ob}}\,{\Cal{C}}$.

Thus it remains to show that also conversely, a $T_2$ uniform space with
covering character at most $\lambda _0$ belongs to 
${\text{\rm{Ob}}}\,{\Cal{C}}$. We will
prove this in three steps: 
\newline
1) for complete metric spaces, with the induced uniformities,
\newline
2) for any metric spaces, with the induced uniformities,
\newline
3) for any $T_2$ uniform spaces.
 
{\bf{3.1.}}
Let $M$ be a complete metric space with cov\,char\,$M \le \lambda _0$.
By Lemma 2 there is a dense,
uniformly continuous map from a closed subspace of a countable product (taken
in ${\bold{T_2Unif}}$) 
$\prod _{n=1}^{\infty } D_{\lambda _n}$ to $M$ --- where 
$D_{\lambda _n}$ is a 
discrete uniform space of cardinality $\lambda _n \,\,(< \lambda _0)$. 

Therefore we have for each $n$ that 
$D_{\lambda _n} \in {\text{\rm{Ob}}}\,{\Cal{C}}$, hence
$\prod _{n=1}^{\infty } D_{\lambda _n} \in $ ${\text{\rm{Ob}}}\,{\Cal{C}}$, 
hence all closed
subspaces of $\prod _{n=1}^{\infty } D_{\lambda _n}$ belong to 
${\text{\rm{Ob}}}\,{\Cal{C}}$, 
as well as all dense images of these closed subspaces 
belong to ${\text{\rm{Ob}}}\,{\Cal{C}}$. 
Therefore $M \in {\text{\rm{Ob}}}\,{\Cal{C}}$ for any
complete metric space $M$ with cov\,char\,$M \le \lambda _0$, with the induced
uniformity.

{\bf{3.2.}}
Let $(M,\varrho )$ be a metric space with cov\,char\,$M \le \lambda _0$. 
As in the proof of Lemma 2 we may assume diam\,$M < 1$.
Then its completion $\gamma (M,
\varrho ) =: (\gamma M, {\tilde{ \varrho }})$ 
has the same covering character, hence, by {\bf{3.1}}, 
belongs to ${\text{\rm{Ob}}}\,{\Cal{C}}$.
Let $m_0 \in \gamma M$ be arbitrary, but fixed. 
Then $(\gamma M) \setminus \{ m_0 \} =
f^{-1}\left( (0,1] \right) $, where $f: \gamma M  \to [0,1]$ is defined as 
$f(m) := {\tilde {\varrho }}(m_0,m)$, 
for each $m \in \gamma M$. Since $\gamma M \in 
{\text{\rm{Ob}}}\,{\Cal{C}}$ and 
$(0,1] \in {\text{\rm{Ob}}}\,{\Cal{C}}$, therefore, 
by 2) of Lemma 1, $(\gamma M) \setminus \{ m_0 \} =f^{-1} \left( (0,1] \right)
\in 
{\text{\rm{Ob}}}\,{\Cal{C}}$. 
Then 1) of Lemma 1 implies that any subspace of $\gamma M$
belongs to $ {\text{\rm{Ob}}}\,{\Cal{C}}$. 
In particular, $M \in {\text{\rm{Ob}}}\,{\Cal{C}}$, for any
metric space $M$ with cov\,char\,$M \le \lambda _0$, with the induced
uniformity.

{\bf{3.3.}}
Let $X$ be a $T_2$ uniform space with cov\,char\,$X \le \lambda _0$. Then $X$
is a subspace of a product of metric spaces $M_{\alpha }$, for $\alpha \in A$.
We may suppose that the restriction of each projection $\pi _{\alpha }: \prod
_{\alpha \in A} M_{\alpha } \to M_{\alpha }$ to $X$ is surjective. Then, for
each $\alpha \in A$, we have cov\,char\,$M_{\alpha } \le \lambda _0$, hence, 
by {\bf{3.2}}, $M_{\alpha } \in {\text{\rm{Ob}}}\,{\Cal{C}}$. Now let us
embed each
$M_{\alpha }$ to $N_{\alpha } := M_{\alpha } \times [0,1]$, via $m_{\alpha }
\mapsto (m_{\alpha }, 0)$, for $m_{\alpha } \in M_{\alpha }$. Then
cov\,char\,$N_{\alpha }=$ cov\,char\,$M_{\alpha } \le \lambda _0$, for each
$\alpha \in A$.

Let $n_{\alpha } \in N_{\alpha }$ be arbitrary. Then cov\,char\,$(N_{\alpha }
\setminus \{ n_{\alpha } \} ) 
= $ cov\,char\,$N_{\alpha } \le \lambda _0$, hence,
by {\bf{3.2}}, 
$$N_{\alpha } \setminus \{ n_{\alpha } \}  \in 
{\text{\rm{Ob}}}\,{\Cal{C}}, {\text{ and therefore }} \prod _{\alpha \in A} 
(N_{\alpha } \setminus \{ n_{\alpha } \} ) \in 
{\text{\rm{Ob}}}\,{\Cal{C}}.
$$ 
Since $n_{\alpha } $ is not an isolated point
of $N_{\alpha }$, therefore $\prod _{\alpha \in A} 
(N_{\alpha } \setminus \{ n_{\alpha } \} )$ is dense in 
$\prod _{\alpha \in A} N_{\alpha }$, hence any subspace of 
$\prod _{\alpha \in A} N_{\alpha }$, containing $\prod _{\alpha \in A} 
(N_{\alpha } \setminus \{ n_{\alpha } \} )$ (as a dense subspace), belongs to 
$ {\text{\rm{Ob}}}\,{\Cal{C}}$.
In particular, 
$$
(\prod _{\alpha \in A} N_{\alpha }) \setminus \{ \langle 
n_{\alpha } \rangle \} \in 
{\text{\rm{Ob}}}\,{\Cal{C}}, {\text{ for arbitrary }} \langle 
n_{\alpha } \rangle \in \prod _{\alpha \in A} N_{\alpha }.
$$ 

\newpage

By 1) of Lemma 1,
then any subspace of $\prod _{\alpha \in A} N_{\alpha }$ belongs to 
$ {\text{\rm{Ob}}}\,{\Cal{C}}$. In particular, 
any subspace of $\prod _{\alpha \in A} M_{\alpha }$, e.g., $X$, belongs to 
$ {\text{\rm{Ob}}}\,{\Cal{C}}$, for any $T_2$ uniform
space $X$ with cov\,char\,$X \le \lambda _0$.
 
Together with the first paragraph of {\bf{3}} this gives that 
${\text{\rm{Ob}}}\,{\Cal{C}}=
\{ X \in {\text{\rm{Ob}}}\,{\bold{T_2Unif}} \mid X$ has a covering
character at most $\lambda _0 \}$.

{\bf{4.}}
There remains the case, from the case distinction in {\bf{2}}, 
when all uniformly discrete spaces $D_{\lambda }$ (of
cardinality $\lambda $) belong to 
$ {\text{\rm{Ob}}}\,{\Cal{C}}$. Then by the above proof,
for any cardinal $\lambda _0$, all $T_2$ uniform spaces 
with covering character
at most $\lambda _0$ belong to 
$ {\text{\rm{Ob}}}\,{\Cal{C}}$. That is, all $T_2$ uniform
spaces belong to 
$ {\text{\rm{Ob}}}\,{\Cal{C}}$, hence ${\Cal{C}}={\bold{T_2Unif}}$.

{\bf{5.}}
There remains to give four examples. These are (except the second one)
the same as in \cite{P}. All but
the 1-st, 2-nd, 3-rd or 
4-th properties are satisfied by the subclasses of ${\bold{T_2Unif}}$
consisting of finite spaces,
of spaces with connected topology, of spaces where the closure of any at most
countably infinite set is compact, or of $0$-dimensional compact $T_2$ uniform
spaces, respectively.
$\blacksquare $
\enddemo


\demo{Proof of Theorem 3}
We will follow the proofs of Theorems 1 and 2.
We only need to prove $1) \Longrightarrow 2)$.

{\bf{1.}}
Like in the proof of Theorem 1, second paragraph, we see that 
$\{ X \in {\text{{\rm{Ob}}}}\,{\bold{Unif}} \mid |X| \le 1 \} \subset 
{\text{\rm{Ob}}}\,{\Cal{C}}$. If here we have equality, we are done.
 
Now suppose that here we do not have equality, i.e., 
$ {\text{\rm{Ob}}}\,{\Cal{C}} $ contains a
uniform space $X$ with at least two points. Then $X$ has a 
subspace consisting of two points, thus some two-point space 
belongs to $ {\text{\rm{Ob}}}\,{\Cal {C}}$. 
Then its bijective image, the two-point indiscrete
space $I_2$ also belongs to $ {\text{\rm{Ob}}}\,{\Cal {C}}$.
Then any subspace of any power of $I_2$ belongs to 
$ {\text{\rm{Ob}}}\,{\Cal{C}}$, hence
$\{ X \in {\text{\rm{Ob}}}\,{\bold{Unif}} \mid X$ is indiscrete$ \} \subset
{\text{\rm{Ob}}}\,{\Cal{C}}$. If here we have equality, we are done.

Now suppose that here we do not have equality, i.e., ${\Cal{C}} $ contains a
non-indiscrete uniform space $X$. 
Then $X$ has a 
subspace consisting of two points, that is a discrete two-point space $D_2$,
and that has to belong 
to $ {\text{\rm{Ob}}}\,{\Cal {C}}$. Then all subspaces of all finite powers
of $D_2$ belong to $ {\text{\rm{Ob}}}\,{\Cal{C}}$, i.e.,  
$\{ X \in {\text{{\rm{Ob}}}}\,{\bold{Unif}} \mid X $ is a finite discrete 
space$ \} \subset {\text{\rm{Ob}}}\,{\Cal{C}}$.

{\bf{2.}}
The class of cardinalities $\lambda $ (finite or infinite), 
for which the discrete space $D_{\lambda
}$ of cardinality $\lambda $ belongs to 
$ {\text{\rm{Ob}}}\,{\Cal{C}}$, forms an initial
segment of all cardinalities, i.e., it is of the form $\{ \lambda \mid \lambda
< \lambda _0 \} $, or $\{ \lambda \mid \lambda $ is a cardinal$\} $. In the
first case, by the last sentence of
{\bf{1}}, we have $\lambda _0 \ge \aleph _0$.

{\bf{3.}}
We begin with the case when this initial segment is $\{ \lambda \mid \lambda
< \lambda _0 \} $. No uniform space in $ {\text{\rm{Ob}}}\,{\Cal{C}}$ 
can have a
covering character greater than $\lambda _0$, since such a space contains a
(closed) subspace $D_{\lambda _0}$, and then we would have $D_{\lambda _0} \in
{\text{Ob}}\,{\Cal{C}}$.

Thus it remains to show that, also conversely, a uniform space with
covering character at most $\lambda _0$ belongs to 
$ {\text{\rm{Ob}}}\,{\Cal{C}}$. 

Let $\lambda < \lambda _0$, and let $I$ be an indiscrete space. Then by
{\bf{1}} and {\bf{2}} we have that 
$I, D_{\lambda } \in {\text{\rm{Ob}}}\,{\Cal{C}}$, hence
each subspace of $D_{\lambda } \times I$
belongs to $ {\text{\rm{Ob}}}\,{\Cal{C}}$. 
That is, each uniform space, with underlying set
$X$, say,
that has a covering
base consisting of a single partition $P$ of cardinality $\lambda $, 
belongs to $ {\text{\rm{Ob}}}\,{\Cal{C}}$. Let us denote this space $X$ 
by $X_P$.

Now let $Y$ be a uniform space with underlying set $X$,
having a covering base ${\Cal{P}}$
consisting of all partitions $P$ of $X$, 
of cardinalities $|P| < \lambda _0$.
Then $Y$ can be embedded to $\prod _{P \in {\Cal P}} X_P\,\,(\in 
{\text{\rm{Ob}}}\,{\Cal{C}}$) 
via the diagonal map. That is, the subspace of this product space, which is the
diagonal, is isomorphic to $Y$. Therefore also $Y \in 
{\text{\rm{Ob}}}\,{\Cal{C}}$.

\newpage

Of course, $Y$ has another covering base, consisting of all covers of $X$ of
cardinalities less than $\lambda _0$. This implies that any uniform structure 
on the underlying set $X$, having a covering base consisting of covers of
cardinalities less than $\lambda _0$, is a bijective image of $Y$, hence 
belongs to $ {\text{\rm{Ob}}}\,{\Cal{C}}$ 
as well. That is, any uniform space, with
covering character at most $\lambda _0$, belongs to 
$ {\text{\rm{Ob}}}\,{\Cal{C}}$. 

Together with the first paragraph of {\bf{3}} this gives that 
${\text{\rm{Ob}}}\,{\Cal{C}}=
\{ X \in {\text{\rm{Ob}}}\,{\bold{Unif}} \mid X$ has a covering
character at most $\lambda _0 \}$.

{\bf{4.}}
There remains the case, from the case distinction in {\bf{2}}, 
when all uniformly discrete spaces $D_{\lambda }$ (of
cardinality $\lambda $) belong to 
${ {\text{\rm{Ob}}}\,\Cal{C}}$. Then by the above proof,
for any cardinal $\lambda _0$, all uniform spaces with covering character
at most $\lambda _0$ belong to $ {\text{\rm{Ob}}}\,{\Cal{C}}$. 
That is, all uniform
spaces belong to 
$ {\text{\rm{Ob}}}\,{\Cal{C}}$, hence ${\Cal{C}}={\bold{Unif}}$.

{\bf{5.}}
There remains to give three examples.
The first one is the same as in \cite{P}. All but
the 1-st, 2-nd or 3-rd 
properties are satisfied by the subclasses of ${\bold{Unif}}$ consisting 
of finite spaces,
of spaces with connected topology, or by $T_2$ uniform spaces,
respectively.
$\blacksquare $
\enddemo


Before the proof of Theorems 4, 5, 6 we give a lemma. Lemma 3, 1) is surely
known (algebraic subcategories of ${\bold{Set}}$), 
but could not locate it, therefore we give its simple proof.

\proclaim{Lemma 3}
Let ${\bold{T}}$ be ${\bold{Top}}$, ${\bold{Prox}}$ or ${\bold{Unif}}$. Let 
${\Cal{C}} \subset {\bold{T}}$ be algebraic. 

1) If \,{\rm{Ob}}\,${\Cal{C}} \subset 
\{ X \in $ {\rm{Ob}}\,${\bold{T}} \mid X$ is indiscrete$\}$, then 
{\rm{Ob}}\,${\Cal{C}} = 
\{ X \in $ {\rm{Ob}}\,${\bold{T}} \mid |X| = 1 \} $, 
or \,{\rm{Ob}}\,${\Cal{C}} =
\{ X \in $ {\rm{Ob}}\,${\bold{T}} \mid |X| \le 1 \} $, 
or \,{\rm{Ob}}\,${\Cal{C}} = 
\{ X \in $ {\rm{Ob}}\,${\bold{T}} \mid X$ is indiscrete$\} $.

2) If \,${\bold{T}}$ is ${\bold{Prox}}$ or ${\bold{Unif}}$, and 
${\Cal{C}}$ is closed under products and closed subspaces, then either
{\rm{Ob}}\,${\Cal{C}} = \{ X \in $ {\rm{Ob}}\,${\bold{T}} \mid X$ 
is indiscrete$\} $, or 
{\rm{Ob}}\,${\Cal{C}} \subset \{ X \in $ {\rm{Ob}}\,${\bold{T}} \mid X$ 
is $T_2 \} $.
\endproclaim


\demo{Proof}
{\bf{1.}}
We begin with the proof of 1). 

The category
${\Cal{C}}$, as a category in its own right, has products, which are preserved
by the underlying set functor. Therefore 
the empty product, the one-point algebra
(space) belongs to $ {\text{\rm{Ob}}}\,{\Cal{C}}$. 
Then either $\emptyset \not\in 
{\text{\rm{Ob}}}\,{\Cal{C}}$ or 
$\emptyset \in {\text{\rm{Ob}}}\,{\Cal{C}}$. Therefore 
$$
\{ X \in {\text{\rm{Ob}}}\,{\bold{T}} \mid |X|=1 \} \subset 
{\text{\rm{Ob}}}\,{\Cal{C}}
\,\,(\not\ni \emptyset ), {\text{ or }} 
\{ X \in {\text{\rm{Ob}}}\,{\bold{T}} \mid |X| \le 1 \} 
\subset {\text{\rm{Ob}}}\,{\Cal{C}}.
$$
If in one of these inclusions we have equality, we are done.
 
Therefore we may suppose that some space $X$ belongs to 
${\text{\rm{Ob}}}\,{\Cal{C}}$, where $|X| \ge 2$.

Then $X$ is indiscrete, and
all powers $X^{\alpha }$ of $X$, 
taken in ${\Cal{C}}$, belong to $ {\text{\rm{Ob}}}\,{\Cal{C}}$. Now, the
underlying set functor $U$ of ${\Cal{C}}$
preserves products, hence these products are the
indiscrete structures on $(UX)^{\alpha }$ (i.e., the powers taken in
${\bold{T}}$), hence
indiscrete spaces of arbitrarily large cardinality belong to $
{\text{\rm{Ob}}}\,{\Cal{C}}$.
Let us consider $X^{\alpha } \in {\text{\rm{Ob}}}\,{\Cal{C}}$. 
Let us consider any 
subset $Y$ of $X^{\alpha }$. 
Let $x_1,x_2 \in X$, with $x_1 \ne x_2$ (recall $|X| \ge 2$).
Let us consider the ${\bold{T}}$-morphisms (hence ${\Cal{C}}$-morphisms)
$f,g: X^{\alpha } \to
X$, defined by $f(z)=x_1$ for all $z \in X^{\alpha }$, and $g(z)=x_1$ for all 
$z \in Y$ and $g(z)=x_2$ for all $z \in X^{\alpha } \setminus Y$. Recall that
the equalizer of $f,g$ is preserved by the underlying set functor of
${\Cal{C}}$, hence $Y \in {\text{\rm{Ob}}}\,{\Cal{C}}$ (up to isomorphy, but
$ {\text{\rm{Ob}}}\,{\Cal{C}}$ is isomorphism closed).
Therefore all indiscrete spaces of cardinality at most
$|X^{\alpha }|$ belong to $ {\text{\rm{Ob}}}\,{\Cal{C}}$.
Hence
$\{ C \in {\text{{\rm{Ob}}}}\,{\bold{T}} \mid C$ is indiscrete$ \} \subset
{\text{\rm{Ob}}}\,{\Cal{C}}$. 
Since the converse inclusion holds by 

\newpage

hypothesis, we have here in
fact equality. 
(Cf. also the proof of \cite{R91a}, Proposition 1.1.)

{\bf{2.}}
We turn to the proof of 2).

If $ {\text{\rm{Ob}}}\,{\Cal{C}} \subset \{ X \in 
{\text{\rm{Ob}}}\,{\bold{T}} \mid X$ is $T_2 \} $, we are done. Therefore
let $ {\text{\rm{Ob}}}\,{\Cal{C}} \not\subset \{ X \in 
{\text{\rm{Ob}}}\,{\bold{T}} \mid X$ is $T_2 \} $,
and let us choose $C \in {\Cal{C}}$ that is not $T_2$, i.e., that is not
$T_0$. Then some point of $C$ is not closed, and its closure, $X$, say, is
a closed indiscrete subset of $C$, with $|X| > 1 $. Then, by closed
hereditariness and productivity of ${\Cal{C}}$, we have 
$\emptyset , X \in {\text{\rm{Ob}}}\,{\Cal{C}}$, and 
also any power $X^{\alpha } $ belongs to 
$ {\text{\rm{Ob}}}\,{\Cal{C}}$. Then repeating 
the considerations in 1) 
we obtain
$\{ X \in {\text{{\rm{Ob}}}}\,{\bold{T}} \mid X$ is indiscrete$ \} \subset 
{\text{\rm{Ob}}}\,{\Cal{C}}$. If here we have equality, we are done.

Therefore we may assume that some non-indiscrete space $C$ belongs to
$ {\text{\rm{Ob}}}\,{\Cal{C}}$. Then the indiscrete space $I$ on $UC$ also
belongs to 
$ {\text{\rm{Ob}}}\,{\Cal{C}}$, and we have 
a bijection $b: C \to I$ that is not an isomorphism 
in ${\bold{T}}$, hence it is not an isomorphism in ${\Cal{C}}$ either. 
However, for
an algebraic category ${\Cal{C}}$, the underlying set functor $U$ reflects
isomorphisms, and we have a contradiction.
$\blacksquare $
\enddemo


\demo{Proof of Theorem 4}
We only need to prove $1) \Longrightarrow 2)$.

{\bf{1.}}
The left adjoint of the underlying set functor
$U:{\Cal{C}} \to {\bold{Set}}$ will be denoted by $F$. 
Objects of ${\bold{T}}$ will be called {\it{spaces}}, and if we
investigate an object of ${\Cal{C}}$, it will be called an {\it{algebra}}.

Like in the proof of Theorem 1, second paragraph, we see that 
$$
\{ X \in {\text{\rm{Ob}}}\,{\bold{T}} \mid |X| \le 1 \} \subset 
{\text{\rm{Ob}}}\,{\Cal{C}}. \tag *
$$ 
If here we have equality, we are done. 
Therefore let 
$ {\text{\rm{Ob}}}\,{\Cal{C}}$ contain an object $C$ with $|UC| \ge 2$. 

{\bf{2.}}
We distinguish two cases: 

1) $ {\text{\rm{Ob}}}\,{\Cal{C}} \subset \{$indiscrete spaces in
$ {\text{\rm{Ob}}}\,{\bold{T}} \} $, or

2) there exists 
$C \in {\text{\rm{Ob}}}\,{\Cal{C}}$, such that $C$ (as an object of 
${\bold{T}}$) is not indiscrete.

{\bf{3.}} 
In the first case, by Lemma 3, 1) and \thetag{*}, we have 
${ {\text{\rm{Ob}}}\,\Cal{C}} = 
\{ X \in {\text{\rm{Ob}}}\,{\bold{T}} \mid |X| \le 1 \}$ or
$ {\text{\rm{Ob}}}\,{\Cal{C}} = \{ X \in $ Ob\,${\bold{T}} \mid X $ is
indiscrete$\}$.

{\bf{4.}}
In the second case 
$ {\text{\rm{Ob}}}\,{\Cal{C}}$ contains a non-indiscrete object $C$,
hence as its subspace, also contains a non-indiscrete object with two points.
Hence we may suppose that $|UC|=2$. Then, for ${\bold{T}}={\bold{Prox}}$
and  ${\bold{T}}={\bold{Unif}}$
we have that $C$ is the two-point discrete space $D_2$. 
For ${\bold{T}}={\bold{Top}}$ we have
that $C$ is the two-point discrete space, or the Sierpi\'nski space. However,
observe that the square of the Sierpi\'nski space contains a two-point
discrete subspace, so we may assume that $C$ is the two-point discrete space
$D_2$ for ${\bold{T}}={\bold{Top}}$ as well.

Let $UC = UD_2 = \{c_1, c_2 \} $. Let $S$ be a set with $|S| \ge 2$,
and let us consider the free algebra $FS \in $ 
Ob\,${\Cal{C}}$. We have the unit of adjunction $\eta _S:S \to UFS$.
Then for any ${\bold{Set}}$-morphism $f:S \to \{ c_1,c_2 \} = UC$, 
there exists a ${\Cal{C}}$-morphism 
$\varphi : FS \to C$ such that $f=(U \varphi ) \circ
\eta _S$. This readily implies that $\eta _S$ is an injection. Moreover it
also implies that
$\eta _S S \,\,(\subset UFS)$, considered as a subspace of $FS$, that (by
hereditariness of ${\Cal{C}}$)
satisfies $\eta _S S \in $ Ob\,${\Cal{C}}$, also satisfies the following. It
has a topology/proximity finer than (thus equal to) the
one projectively generated by all ${\bold{Set}}$-morphisms 
to $\{ c_1,c_2 \} = UC$, i.e., 
the discrete topology/proximity on $U(\eta _S S)$. Or it has 
a uniformity finer than
the finest precompact uniformity on $U( \eta _S S)$ (i.e., the one having as
a covering base all finite partitions of $U( \eta _S S)$), respectively.

\newpage

Now let ${\bold{T}}={\bold{Top}}$ or ${\bold{T}}={\bold{Prox}}$. Then,
by epireflectivity, this discrete subspace $\eta_S S$ belongs to 
Ob\,${\Cal{C}}$ (and the discrete spaces of cardinality at most $1$ belong to
Ob\,${\Cal{C}}$ by \thetag{*}).
Hence
Ob\,${\Cal{C}} \supset \{$discrete spaces in ${\bold{T}} \} $.
For ${\bold{T}}={\bold{Unif}}$ the same reasoning 
gives only that the uniformity $\eta _S S$
on $U(\eta _S S)$, finer than the finest
precompact uniformity on $U(\eta _S S)$, belongs to Ob\,${\Cal{C}}$.

The two-point discrete space $D_2$ in ${\bold{T}}$ belongs to 
Ob\,${\Cal{C}}$. Hence, 
by epireflectivity, also $D_2^{\aleph _0}
\in $ Ob\,${\Cal{C}}$, where the power is taken in ${\bold{T}}$. 
However, $D_2^{\aleph _0}$ is the Cantor set, or the Cantor set with
its unique compatible proximity, or the Cantor set with its unique compatible
uniformity, respectively.
Let 
$$
S:=U(D_2^{\aleph _0})= (UD_2)^{\aleph _0}.
$$ 
For simplicity, we assume that the embedding $\eta _S:S \to UFS$ is pointwise
identical. Then for ${\bold{T}}={\bold{Top}}$ and 
${\bold{T}}={\bold{Prox}}$ the space $\eta _S S$ (as a subspace of $FS$)
is discrete, hence is
strictly finer than
$D_2^{\aleph _0}$. 
For ${\bold{T}}={\bold{Unif}}$ the space
$\eta _S S$ is finer than the finest precompact uniformity on $U(\eta _S S)$.
In all three cases the space $\eta _S S$ is strictly finer than $D_2^{\aleph
_0}$. Thus the identical bijection $b: \eta _S S \to D_2^{\aleph _0}$ is not
an isomorphism in ${\bold{T}}$, hence not an isomorphism in ${\Cal{C}}$
either. However, for
an algebraic category ${\Cal{C}}$, the underlying set functor $U$ reflects
isomorphisms, and we have a contradiction.
Hence case 2) in {\bf{2}}
from our case distinction cannot exist.

{\bf{5.}}
There remains to give three examples. 
All but the 1-st, 2-nd or 3-rd
properties are
satisfied by examples 1) and 2) of \cite{HS71} and by 3) the category of
$T_2$ topological,
proximity or uniform spaces 
(observe that their underlying set functors
do not reflect isomorphisms therefore they are not algebraic).
In 1) we mean discrete topological,
proximity or uniform spaces, which form even a varietal category \cite{HS71}.
In 2) we mean powers of a compact $T_2$
topological space, or of the same space with its unique compatible
proximity or uniformity, consisting of more than one point, that is strongly
rigid --- i.e., whose only continuous self-maps are the identity and the
constant maps (such spaces exist, cf. \cite{HS71}) --- 
together with the empty space, which category is even varietal,
cf. \cite{R92a}, p. 368.
$\blacksquare $
\enddemo


A large part of the next proof is taken from \cite{R82} and \cite{R85a}.

\demo{Proof of Theorem 5}
{\bf{1.}}
The implication $2) \Rightarrow 1)$ is evident in the first case.
For the second case observe that the category of compact $T_2$ uniform spaces 
is canonically
concretely isomorphic to the category of compact $T_2$ topological
spaces (via induced topology/unique compatible uniformity). 
Thus our implication reduces to the analogous implication for the 
category of compact $T_2$ topological spaces, that follows from \cite{R82},
Corollary 1.5.

{\bf{2}}.
We repeat {\bf{1}}, {\bf{2}}, {\bf{3}} from the proof of Theorem 4 word for
word. However, observe that $\{ X \in $ Ob\,${\bold{Unif}} \mid |X|
\le 1 \} $ is an epireflective subcategory of 
$\{ X \in $ Ob\,${\bold{Unif}} \mid X {\text{ is compact }} T_2 \} $.

By the first sentence of {\bf{4}} from the proof of Theorem 4,
Ob\,${\Cal{C}}$ contains a non-indiscrete object $C$.

Then by Lemma 3, 2) we have  
Ob\,${\Cal{C}} \subset {\text{\rm{Ob}}}\,{\bold{T_2Unif}}$.

{\bf{3}}.
From now on we follow \cite{R82}. 

Lemma 1.1 of \cite{R82} will become the following. 
Let every bijection in ${\Cal{C}}$ be a uniform
isomorphism, let $\beta $ be a limit ordinal, and let 
$[0, \beta ]$ be the usual (compact)
ordinal space, 
with the unique uniformity compatible with its order topology. 
Further, let ${\Cal{U}}$ be 
a uniformity on the ordinal
space $[0,\beta )$, which is finer than the 

\newpage

precompact uniformity inherited
from its compactification $[0,\beta ]$ (i.e., the coarsest --- precompact --- 
uniformity compatible with its order topology). Then 
$$ 
[0, \beta ] \in {\text{\rm{Ob}}}\,{\Cal{C}} \Longrightarrow 
([0, \beta ),{\Cal{U}}) \not\in {\text{\rm{Ob}}}\,{\Cal{C}}\,.
$$
The proof remains the same.

In the statement of Proposition 1.2 of \cite{R82}, 
${\bold{Top}}$ has to be replaced by ${\bold{Unif}}$, and of
course, $D_2$ is now a discrete uniform space on two points. 
The proof remains the same, of course replacing topological products and
coproducts by uniform ones.

The statement of Lemma 1.3 of \cite{R82} remains word for word the same, 
of course replacing ${\bold{Top}}$ by ${\bold{Unif}}$, and also 
its proof remains the same.

The assertion of Theorem 1.4 of \cite{R82}
remains the same, of course replacing ${\bold{Top}}$ by ${\bold{Unif}}$.
In the proof the following
changes have to be made. Everywhere, rather than the ordinal space 
$[0, \alpha ]$, with its order topology, we consider 
the unique uniformity compatible with its order topology. 
Moreover, rather than the ordinal spaces 
$[0, \beta )$ or $[\alpha +1, \beta )$, with their order topologies 
(where $\alpha < \beta $), 
we consider the respective compatible uniformities on them,
that are the above mentioned precompact uniformities inherited from their
compactifications $[0, \beta ]$ or $[\alpha +1, \beta ]$ (i.e., their
coarsest compatible uniformities).
For $\alpha < \beta $, the set $[0, \alpha ]$ is not just clopen in $[0,
\beta )$, but together with its complement form a uniform cover of $[0,
\beta )$. Accordingly, topological coproduct at this place 
is replaced by uniform coproduct.

Then the statements of Corollaries 1.5 and 1.6 of \cite{R82}
remain word for word the same,
and also their proofs carry over. (Actually, for \cite{R82}, 
Corollary 1.5, after
the first step of its proof, namely that Ob\,${\Cal{C}} \subset \{ X \in $
Ob\,${\bold{T_2Unif}} \mid X$ is compact $(T_2) \} $,
we can use the canonical concrete
isomorphism of the categories of compact $T_2$ topological and compact $T_2$
uniform spaces mentioned in {\bf{1}} of this proof, 
and then just we have to apply the result of \cite{R82},
Corollary 1.5,
not repeat its proof.)

{\bf{4.}}
There remains to give three examples (in ${\bold{Unif}}$). 
These are the same as in {\bf{5}} of the
proof of Theorem 4, of course in the third case meaning only the category
${\bold{T_2Unif}}$. 
Observe that examples 1) and 2) are even varietal, and
examples 1) (\cite{HS71}) and 3) are closed even for any subspaces.
$\blacksquare $
\enddemo


A large part of the next proof is taken from \cite{HS71} and \cite{R85a}.

\demo{First proof of Theorem 6}
We only need to prove $1) \Longrightarrow 2)$.

{\bf{1.}}
Like in the proof of Theorem 1, second paragraph, we see that 
$\{ X \in {\text{{\rm{Ob}}}}\,{\bold{Unif}} \mid |X| \le 1 \} \subset $
Ob\,${\Cal{C}}$. If here we have equality, we are done.

Therefore we may suppose that some space $X$ belongs to Ob\,${\Cal{C}}$,
where $|X| \ge 2$.

{\bf{2.}}
Now we make a case distinction. Either 

1) Ob\,${\Cal{C}} \not\subset $ Ob\,${\bold{T_2Unif}}$, or

2) Ob\,${\Cal{C}} \subset $ Ob\,${\bold{T_2Unif}}$.

{\bf{3.}}
We begin with case 1). Then, by Lemma 3, 2), we have
Ob\,${\Cal{C}}=
\{ X \in {\text{{\rm{Ob}}}}\,{\bold{Unif}} \mid X$ is indiscrete$ \} $.
(Cf. also the proof of \cite{R91a}, Proposition 1.1.)

{\bf{4.}}
We turn to case 2), i.e., when 
Ob\,${\Cal{C}} \subset $ Ob\,${\bold{T_2Unif}}$.
Then we can repeat the proof of Lemmas 1 and 2, Corollaries 1 and 2 and
the Theorem from \cite{HS71}. 

\newpage

We only have to change the words topological spaces,
continuous maps, topological quotient spaces and maps etc. to their uniform
counterparts.
Thus we obtain that 
$$
{\text{Ob}}\,{\Cal{C}}=
\{ X \in {\text{{\rm{Ob}}}}\,{\bold{Unif}} \mid |X| \le 1\} , {\text{ or }}
{\text{Ob}}\,{\Cal{C}}=
\{ X \in {\text{{\rm{Ob}}}}\,{\bold{Unif}} \mid X {\text{ is compact }} T_2  
\} . 
$$

{\bf{5.}}
There remains to give three examples (in ${\bold{Unif}}$). 
These are the same as in {\bf{5}} of the
proof of Theorem 4, of course in the third case meaning only the category
${\bold{T_2Unif}}$.
Observe that
examples 1) (\cite{HS71}) and 3) are closed even for any subspaces.
$\blacksquare $
\enddemo


Of course, Theorem 6 also
follows from Theorem 5. However, this proof of Theorem 6,
although is shorter to write, is in fact more complicated. Namely, in the above
first proof of Theorem 6 we used the proof from \cite{HS71}, 
which is simpler than the 
proof from \cite{R82} of a more general theorem,
used in the proof of Theorem 5.


\demo{Second proof of Theorem 6}
We only deal with $1) \Longrightarrow 2)$.

Parts {\bf{1}}, {\bf{2}}, {\bf{3}} from the above proof of Theorem 6 are just
copied. Thus, in particular, we assume that some uniform space $X$ belongs to 
Ob\,${\Cal{C}}$, where $|X| \ge 2$, and Ob\,${\Cal{C}} \subset
{\text{Ob}}\,{\bold{T_2Unif}}$. 

Since a varietal category is algebraic, we have by
Theorem 5 that ${\Cal{C}}$ is an epireflective subcategory of compact $T_2$
uniform spaces. By the canonical concrete isomorphism of the categories
of compact $T_2$
uniform spaces and compact $T_2$ topological spaces (from {\bf{1}} of the
proof of Theorem 5), we have a corresponding 
epireflective and varietal
subcategory ${\Cal{C}}'$ of compact $T_2$ topological spaces.
Then \cite{HS71} Theorem or
\cite{R82}, Corollary 1.6 implies that 
$$
\cases
{\text{Ob}}\,{\Cal{C}}'= \{
{\text{compact }} T_2 {\text{ topological spaces}} \} , 
\\ 
{\text{ hence }}  
{\text{Ob}}\,{\Cal{C}}= \{
{\text{compact }} T_2 {\text{ uniform spaces}} \} .
\endcases
$$
$\blacksquare $
\enddemo


{\bf{Acknowledgements.}} The author expresses his thanks to D. Petz, for useful
conversations on the subject of the paper, and to G. Richter, for kindly
having sent him some of his papers.

\enddocument